# On Disturbance-to-State Adaptive Stabilization without Parameter Bound by Nonlinear Feedback of Delayed State and Input

**Iasson Karafyllis[*], Miroslav Krstic[**] and Alexandros Aslanidis[*]**

[*]Dept. of Mathematics, National Technical University of Athens,
Zografou Campus, 15780, Athens, Greece,
email: iasonkar@central.ntua.gr; iasonkaraf@gmail.com

[**]Dept. of Mechanical and Aerospace Eng., University of California, San Diego, La Jolla, CA 92093-0411, U.S.A., email: krstic@ucsd.edu

## Abstract

We complete the first step towards the resolution of several decades-old challenges in disturbance-robust adaptive control. For a scalar linear system with an unknown parameter for which no a priori bound is given, with a disturbance that is of unlimited magnitude and possibly persistent (not square integrable), and without a persistency of excitation necessarily verified by the state, we consider the problems of (practical) gain assignment relative to the disturbance. We provide a solution to these heretofore unsolved feedback design problems with the aid of infinite-dimensional nonlinear feedback employing distributed delay of the state and input itself. Specifically, in addition to (0) the global boundedness of the infinite-dimensional state of the closed-loop system when the disturbance is present, we establish (1) practical input-to-output stability with assignable asymptotic gain from the disturbance to the plant state; (2) assignable exponential convergence rate; and (3) assignable radius of the residual set. The accompanying identifier in our adaptive controller guarantees (4) boundedness of the parameter estimate even when disturbances are present; (5) an ultimate estimation error which is proportional to the magnitude of the disturbance with assignable gain when there exists sufficient excitation of the state; and (6) exact parameter estimation in finite-time when the disturbance is absent and there is sufficient excitation. Among our results, one reveals a tradeoff between "learning capacity" and "disturbance robustness:" the less sensitive the identifier is to the disturbance, the less likely it is to learn the parameter.

## 1. Introduction

The use of delays in feedback control of finite-dimensional, continuous-time control systems is an idea that has been proposed many years ago in control theory. The idea appears in the book [23] and has been exploited by many researchers: for instance, in [2,10,11] it is shown that the intentional use of delays to a delay-free system can guarantee properties that cannot be guaranteed by delay-free designs. More recently, in [14,15] it was shown that the use of delays in adaptive control can ensure features that cannot be guaranteed by delay-free adaptive control.

However, the adaptive control schemes in [14,15] contain delays but are of a hybrid nature (note the event-triggered character of the adaptive controllers in [14,15]). Thus, the question of robustness with respect to modeling errors arises for the resulting closed-loop systems. Simulations in [14,15] showed strong robustness with respect to various uncertainties, but our efforts to provide theoretical guarantees did not succeed. This, of course, is consistent with the difficulty of achieving robustness in nonlinear closed loops under adaptive control and the inherent limitations hidden therein: see the relevant discussions in [8,9,17,20,22,28] as well as the limitations for robustness that appear in discrete-time systems under adaptive control in [3] (but see also the review paper [4] and references therein). In the supervisory framework of direct adaptive control, practical stability in the presence of external disturbances for plants belonging to known bounded sets has been studied in [6,7,21]. In the present paper no known bound on the plant parameter is assumed.

This is the first paper that proposes the intentional use of delays in a conventional (non-hybrid) fashion for adaptive control. The proposed adaptive feedback law makes the closed loop a Neutral Delay system. Neutral functional differential equations have been studied extensively [1,5,13,18,19,23,25,27]. We study the robustness properties of the closed-loop system when modeling errors are present and make neither assumptions of a priori bounds for the unknown parameters (so that no projection scheme can be applied), nor of a priori bounds on the modeling error (so that no dead zone modification can be applied). Due to the difficulty of the performed study and the novel character of the proposed feedback, as in [3,4], we focus on a scalar system. However, even for the scalar linear case the required analysis is quite demanding: the closed-loop system is a neutral delay system not in Hale's form and standard results for neutral delay systems in Hale's form (see [5,13,25,27]) cannot be used. The reader should notice the difference with the recent work [24]: in [24] the authors apply robust adaptive control schemes (leakage) to a system with discrete delays while here we have a delay-free open-loop system where we apply a novel adaptive scheme with delays.

The present work shows that the intentional use of delays in adaptive control can guarantee features that cannot be guaranteed by existing control schemes. We show practical Input-to-Output Stability (IOS; see [29,12]) for the closed-loop system with respect to the disturbance with assignable asymptotic gain, assignable exponential convergence rate and assignable radius of the residual set. Moreover, we guarantee that the state of the closed-loop system remains bounded when the disturbance is present, in addition to having the property of uniform practical global asymptotic stability with exponential and assignable convergence rate in the absence of disturbance. Let us notice here the difference between the proposed adaptive controller and other conventional nonlinear adaptive controllers proposed in the literature. First, some adaptive control schemes do not guarantee that the state of the closed-loop system is bounded when a disturbance is present. Second, some adaptive schemes guarantee that the state of the closed-loop is bounded when a disturbance of unknown bound is present [8,9,22,28] (e.g., the $\sigma-$modification and other leakage schemes), but the radius of the residual set for the output is neither known nor assignable.

The adaptive feedback law is accompanied by a hybrid identifier, which only provides estimates not used for feedback (and consequently the identifier does not affect the closed-loop system). The identifier guarantees: (i) exact parameter estimation in finite-time when the disturbance is absent and there is sufficient excitation, (ii) boundedness of the parameter estimate even when disturbances are present, and (iii) an ultimate estimation error which is proportional to the magnitude of the disturbance with an assignable gain when there exists sufficient excitation of the state. These features cannot be provided by other identifiers in the literature.

The evident novelty of the proposed adaptive feedback law is accompanied by a novel stability analysis of the closed-loop system which combines Lyapunov arguments and a small-gain analysis (see the proof of Theorem 2 below) in Sobolev spaces. Since we could not find existence/uniqueness results for neutral delay systems not in Hale's form with measurable disturbances, we also prove existence/uniqueness (see Theorem 1 below).

The structure of the paper is as follows. All the main results are stated and discussed in Section 2. Section 3 and Section 4 of the paper are devoted to the proofs of all main results and many auxiliary results. The concluding remarks of the present work are provided in Section 5.

**Notation.** Throughout this paper, we adopt the following notation.

* $\mathbb{R}_+ := [0,+\infty)$. $\mathbb{Z}$ is the set of integers and $\mathbb{N}$ is the set of natural numbers (including zero). For $x \in \mathbb{R}$, $\lceil x \rceil$ denotes the smallest integer that is greater or equal to $x$, i.e., $\lceil x \rceil = \min\{s \in \mathbb{Z} : s \geq x\}$. We use the notation $(x)^+$ for the positive part of the real number $x \in \mathbb{R}$, i.e., $(x)^+ = \max(x,0)$.

* Let $D \subseteq \mathbb{R}^n$ be an open set and let $S \subseteq \mathbb{R}^n$ be a set that satisfies $D \subseteq S \subseteq cl(D)$, where $cl(D)$ is the closure of $D$. By $C^0(S\,;\Omega)$, we denote the class of continuous functions on $S$, which take values in $\Omega \subseteq \mathbb{R}^m$. By $C^k(S\,;\Omega)$, where $k \geq 1$ is an integer, we denote the class of functions on $S \subseteq \mathbb{R}^n$, which take values in $\Omega \subseteq \mathbb{R}^m$ and have continuous derivatives of order $k$. In other words, the functions of class $C^k(S\,;\Omega)$ are the functions which have continuous derivatives of order $k$ in $D = \mathrm{int}(S)$ that can be continued continuously to all points in $\partial D \cap S$. When $\Omega = \mathbb{R}$ then we write $C^0(S)$ or $C^k(S)$.

* Let $I \subseteq \mathbb{R}$ be an interval with non-empty interior. For every $a \geq 1$, $L^a(I)$ denotes the set of equivalence classes of Lebesgue measurable functions $f : I \to \mathbb{R}$ for which $\|f\|_a = \left(\int_I |f(x)|^a\, dx\right)^{1/a} < +\infty$. $L^\infty(I)$ denotes the set of equivalence classes of Lebesgue measurable functions $f : I \to \mathbb{R}$ for which $\|f\|_\infty := \underset{x\in I}{ess\sup}\left(|f(x)|\right) < +\infty$. The scalar product in $L^2(I)$ is denoted by $\langle f, g \rangle = \int_I f(x)g(x)dx$. For every $a \geq 1$, $W^{1,a}(I)$ denotes the set of absolutely continuous functions $f : I \to \mathbb{R}$ with $f, f' \in L^a(I)$. We use the notation $H^1(I)$ for the space $W^{1,2}(I)$.

* Let $r > 0$, $a < b \leq +\infty$ and $u : [a-r, b) \to \mathbb{R}^n$ be given. We use the notation $u_t$ to denote the history at certain $t \in [a,b)$, i.e., $(u_t)(s) = u(t+s)$ for all $s \in [-r, 0]$.

## 2. Main Results

Consider the scalar system

$$\begin{aligned} &\dot{x}(t) = \theta x(t) + u(t) + d(t) \\ &x, u, d \in \mathbb{R} \end{aligned} \tag{2.1}$$

where $\theta \in \mathbb{R}$ is a constant unknown parameter. We assume that $d \in L^{\infty}\left(\mathbb{R}_{+}\right)$ is a disturbance that is not measured. It is evident from (2.1) that for every $r > 0$ the following identity is valid for $t \geq r$ :

$$x^2(t) - x^2(t-r) - 2\left\langle x_t, u_t \right\rangle = 2\theta \left\| x_t \right\|_2^2 + 2\left\langle x_t, d_t \right\rangle$$

Here and in what follows all functions $x_t, u_t, d_t$ are considered to be defined on $[-r, 0]$ and consequently, we have $\left\langle x_t, d_t \right\rangle = \int_{-r}^{0} x(t+s)d(t+s)ds$, $\left\langle x_t, u_t \right\rangle = \int_{-r}^{0} x(t+s)u(t+s)ds$ and $\left\| x_t \right\|_2^2 = \int_{-r}^{0} x^2(t+s)ds$. It is tempting to design a feedback law for (2.1) based on the following estimate for the unknown parameter $\theta \in \mathbb{R}$ when $x_t \neq 0$:

$$\hat{\theta}(t) = \frac{x^2(t) - x^2(t-r) - 2\left\langle x_t, u_t \right\rangle}{2\left\| x_t \right\|_2^2} \tag{2.2}$$

For example, one could think of a feedback law of the form $u(t) = -\left(2c + \hat{\theta}(t)\right)x(t)$, where the parameter estimate $\hat{\theta}(t)$ is given by (2.2) and $c > 0$ is a constant. There are two major issues for the implementation of this idea. First, the feedback law will give the current value of $u(t)$ based on the estimate $\hat{\theta}(t)$ which also depends on the history of $u(t)$. Therefore, the closed-loop system will not be a conventional delay system. Secondly, the mapping that gives the estimate $\hat{\theta}(t)$ is not continuous (it is not even locally bounded). Consequently, there will be severe problems for establishing the well-posedness of the closed-loop system.

We seek a modification of the above idea, which can guarantee the well-posedness of the closed-loop system as well as a set of properties that cannot be guaranteed by other adaptive controllers.

We next investigate the properties of the closed-loop system under the adaptive controller

$$u(t) = -\left(2c + p\left(x_t, u_t\right)\right)x(t) \tag{2.3}$$

$$\frac{d\,\hat{\theta}}{d\,t}(t) = 0 \text{, for all } t \in \left[ir, (i+1)r\right),\ i \in \mathbb{N} \tag{2.4}$$

$$\hat{\theta}((i+1)r) = \hat{\theta}(ir) \text{, if } \max_{ir \leq s \leq (i+1)r}\left(\left\| x_s \right\|_2\right) < \sigma,\ i \in \mathbb{N} \tag{2.5}$$

$$\hat{\theta}((i+1)r) = q\left(x_{\tau},u_{\tau}\right), \text{ if } \max_{ir \le s \le (i+1)r}\left(\|x_s\|_2\right) \ge \sigma,\ i \in \mathbb{N}$$

$$\text{with } \tau = \max\left\{t \in \left[ir,(i+1)r\right] : \|x_t\|_2 = \max_{ir \le s \le (i+1)r}\left(\|x_s\|_2\right)\right\} \tag{2.6}$$

where $\varepsilon, c, r, \sigma > 0$ are constants and

$$p(x,u) := \frac{\left(x^2(0) - x^2(-r) - 2\langle x,u\rangle\right)^+ + crx^2(0)}{2\left(\|x\|_2^2 + \left(\varepsilon^2 - x^2(0)\right)^+\right)}, \text{ for } (x,u) \in C^0\left([-r,0];\mathbb{R}^2\right) \tag{2.7}$$

$$q(x,u) = \frac{x^2(0) - x^2(-r) - 2\langle x,u\rangle}{2\|x\|_2^2}, \text{ for } (x,u) \in C^0\left([-r,0];\mathbb{R}^2\right),\ x \ne 0 \tag{2.8}$$

Some explanations for the feedback law (2.3) and the identifier (2.4)-(2.6) are needed here.

- Notice that the mapping $q$ defined by (2.8) is simply the expression appearing in the right-hand side of (2.2). So, the identifier in (2.6) uses the estimate (2.2) at the most recent time in the interval $\left[ir,(i+1)r\right]$ for which the $L^2\left((-r,0)\right)$ norm is greater or equal to $\sigma$. When there is no time in the interval $\left[ir,(i+1)r\right]$ for which the $L^2\left((-r,0)\right)$ norm is greater or equal to $\sigma$, then the identifier in (2.5) does not update the parameter estimate $\hat{\theta}$.

- It should be noticed that the feedback law (2.3) does not use the parameter estimate $\hat{\theta}$. Therefore, the identifier does not affect the closed-loop system (2.1), (2.3). Instead, using definition (2.7), the feedback law (2.3) can be written as $u(t) = -\left(2c + \frac{\left(x^2(t) - x^2(t-r) - 2\langle x_t,u_t\rangle\right)^+}{2\left(\|x_t\|_2^2 + \left(\varepsilon^2 - x^2(t)\right)^+\right)}\right)x(t) - \frac{crx^3(t)}{2\left(\|x_t\|_2^2 + \left(\varepsilon^2 - x^2(t)\right)^+\right)}$, where it becomes clear that the parameter estimate has been replaced by $\frac{\left(x^2(t) - x^2(t-r) - 2\langle x_t,u_t\rangle\right)^+}{2\left(\|x_t\|_2^2 + \left(\varepsilon^2 - x^2(t)\right)^+\right)}$ (notice the difference with the estimate $\frac{x^2(t) - x^2(t-r) - 2\langle x_t,u_t\rangle}{2\|x_t\|_2^2}$ in (2.2)). On the other hand, the term $-\frac{crx^3(t)}{2\left(\|x_t\|_2^2 + \left(\varepsilon^2 - x^2(t)\right)^+\right)}$ is a (delay-dependent) nonlinear damping term.

Even with this modification, the right-hand side of (2.3) is not bounded when $\left(x_t,u_t\right)$ is in a bounded subset of $C^0\left([-r,0];\mathbb{R}^2\right)$. This can be shown by the sequence of functions

$$x_n(s)=\begin{cases} 0 & s\in\left[-r,\frac{-r}{n+1}\right] \\ \varepsilon\left(\frac{n+1}{r}s+1\right) & s\in\left(\frac{-r}{n+1},0\right] \end{cases}$$, $u_n=0$ with $n\geq 1$, for which it holds that $\|x_n\|_\infty=\varepsilon$, $\|u_n\|_\infty=0$ and $p(x_n,u_n)=3(n+1)\frac{1+cr}{2r}\to+\infty$ as $n\to+\infty$. However, it can be shown (see Lemma 1 in Section 3) that the right-hand side of (2.3) is indeed bounded when $(x_t,u_t)$ is in a bounded subset of $W^{1,\infty}((-r,0))\times C^0([-r,0])$ or a bounded subset of $H^1((-r,0))\times C^0([-r,0])$. Therefore, the closed-loop system (2.1), (2.3) is a system of Functional Differential Equations with state $(x_t,u_t)\in X$, where

$$\begin{aligned} &X:=\left\{\,(x,u)\in S: u(0)=-\left(p(x,u)+2c\right)x(0)\,\right\} \\ &S:=W^{1,\infty}\left((-r,0)\right)\times C^0\left([-r,0]\right) \end{aligned} \tag{2.9}$$

The well-posedness of the closed-loop system (2.1), (2.3) cannot be guaranteed by using standard results in the literature. For example, one can think of obtaining the solution of the closed-loop system (2.1), (2.3) with initial condition $(x_0,u_0)\in X$ as a solution of the delay system

$$\begin{aligned} &\dot{x}(t)=\theta x(t)+u(t)+d(t) \\ &\dot{\xi}(t)=2x(t)u(t) \\ &u(t)=-\left(\frac{\left(x^2(t)-x^2(t-r)-\xi(t)+\xi(t-r)\right)^+ + crx^2(t)}{2\left(\|x_t\|_2^2+\left(\varepsilon^2-x^2(t)\right)^+\right)}+2c\right)x(t) \end{aligned} \tag{2.10}$$

with initial condition $(x_0,\xi_0)\in\left(W^{1,\infty}((-r,0))\right)^2$ that satisfies $\xi(s)=\xi(-r)+2\int_{-r}^{s}x(l)u(l)dl$ for all $s\in[-r,0]$ (notice that $\xi(t)=\xi(-r)+2\int_{-r}^{t}x(l)u(l)dl$ for all $t\geq -r$). Systems like (2.10) with state space being a Sobolev space are studied in [18,19]. In [18,19] systems like (2.10) are characterized as Neutral Functional Differential Equations. However, the use of the existence/uniqueness results in [18,19] would require the assumption that $d\in C^0(\mathbb{R}_+)$ while we would like to study the general case where $d\in L^\infty(\mathbb{R}_+)$. Moreover, system (2.10) is a neutral delay system not in Hale's form and standard results for neutral delay systems in Hale's form (see [5,13,25,27]) cannot be used for its study. Equation (2.3) by itself is a Functional Difference system (see [26]) but its right-hand side is not independent at 0. Therefore, we are forced to prove the following existence/uniqueness result for the closed-loop system (2.1), (2.3).

**Theorem 1:** *Given $\theta\in\mathbb{R}$, $(x_0,u_0)\in X$ and $d\in L^\infty(\mathbb{R}_+)$ there exist $t_{\max}\in(0,+\infty]$ and a unique pair of continuous mappings $x,u:[-r,t_{\max})\to\mathbb{R}$ with $x\in W^{1,\infty}((-r,T))$ for all $T<t_{\max}$, $x(s)=x_0(s)$, $u(s)=u_0(s)$ for all $s\in[-r,0]$ such that (2.1) holds for $t\in[0,t_{\max})$ a.e. and (2.3) holds for all $t\in[0,t_{\max})$. Moreover, if $t_{\max}<+\infty$ then $\limsup_{t\to t_{\max}^-}\left(|u(t)|\right)=+\infty$.*

Based on Theorem 1, we prove the following stabilization result.

**Theorem 2:** *Let $\omega \in (0,c]$ be an (arbitrary) constant with $\exp(\omega r) < 2$. For every $\theta \in \mathbb{R}$ there exist a constant $M > 0$, a function $a \in K_\infty$ and a non-decreasing function $\gamma : \mathbb{R}_+ \to \mathbb{R}_+$ such that for every $(x_0, u_0) \in X$, $d \in L^\infty(\mathbb{R}_+)$ there exists a unique pair of continuous mappings $x, u : \mathbb{R}_+ \to \mathbb{R}$ with $x \in W^{1,\infty}(\mathbb{R}_+)$, $x(s) = x_0(s)$, $u(s) = u_0(s)$ for all $s \in [-r, 0]$ for which (2.1) holds for $t \geq 0$ a.e., (2.3) holds for all $t \geq 0$ and the following estimates hold for all $t \geq 0$*

$$|x(t)| \leq M \exp\left(-c(t-r)^+\right)|x_0(0)| + \varepsilon + \left(1 + \sqrt{r} M \exp\left(-c(t-r)^+\right)\right) \frac{\|d\|_\infty}{2^{1/4} c} \tag{2.11}$$

$$\|u_t\|_\infty \leq \exp(-\omega t) a\left(\|x_0\|_\infty + \|\dot{x}_0\|_\infty + \|u_0\|_\infty\right) + \gamma\left(\|x_0\|_\infty + \|\dot{x}_0\|_\infty + \|u_0\|_\infty + \|d\|_\infty\right)\|d\|_\infty + \varepsilon \rho \tag{2.12}$$

*where* $\rho := 2|\theta| + \dfrac{1}{2\sqrt{2r}} + c(r+3)$.

The stability estimate (2.11) guarantees Global Practical IOS for the output $Y = x$ with respect to the disturbance $d$ with assignable asymptotic gain (equal to $2^{-1/4} c^{-1}$), assignable exponential convergence rate (equal to $c$) and assignable radius of the residual set (equal to $\varepsilon$). Let us notice here the difference between the adaptive controller (2.3) and other conventional nonlinear adaptive controllers proposed in the literature: some adaptive control schemes do not guarantee that the state of the closed-loop system is bounded when the disturbance $d \in L^\infty(\mathbb{R}_+)$ is present and require additional assumptions for the disturbance in order to guarantee boundedness (e.g., $d \in L^2(\mathbb{R}_+) \cap L^\infty(\mathbb{R}_+)$) which cannot be verified in practice. There are some adaptive schemes that can guarantee that the state of the closed-loop is bounded when the disturbance $d \in L^\infty(\mathbb{R}_+)$ is present: see [8,9] (e.g., the $\sigma-$modification and other leakage schemes). However, in such adaptive control schemes the radius of the residual set for the output $Y = x$ is proportional to $|\theta|$ and therefore, is not known (and not assignable). Both estimates (2.11), (2.12) guarantee uniform practical global asymptotic stability with exponential and assignable convergence rate in the disturbance-free case.

The proof of Theorem 2 is based on a novel stability analysis of the closed-loop system which combines Lyapunov arguments and a small-gain analysis (see Section 4 below) in Sobolev spaces.

An additional difference between the adaptive controller (2.3)-(2.6) and other conventional nonlinear adaptive controllers proposed in the literature arises when the identifier (2.4), (2.5), (2.6) is taken into account. This is shown by the following theorem.

**Theorem 3:** *For every $\theta \in \mathbb{R}$, $(x_0, u_0, \hat{\theta}_0) \in X \times \mathbb{R}$, $d \in L^\infty(\mathbb{R}_+)$, there exists a unique triplet of mappings $x, u \in C^0(\mathbb{R}_+)$, $\hat{\theta} \in L^\infty(\mathbb{R}_+)$, with $\hat{\theta}(0) = \hat{\theta}_0$, $x \in W^{1,\infty}(\mathbb{R}_+)$, $x(s) = x_0(s)$, $u(s) = u_0(s)$ for all $s \in [-r, 0]$ for which (2.1) holds for $t \geq 0$ a.e., (2.3) holds for all $t \geq 0$, (2.4) holds for all $t \in [ir, (i+1)r)$, $i \in \mathbb{N}$ and (2.5), (2.6) hold for all $i \in \mathbb{N}$. Moreover, if there exists $t \geq r$ with $\|x_t\|_2 \geq \sigma$ then there exists a finite time $T \geq r$ for which the following estimate holds:*

$$\left|\hat{\theta}(t) - \theta\right| \leq \frac{\sqrt{r}}{\sigma} \|d\|_\infty \text{ for all } t \geq T \tag{2.13}$$

Estimate (2.4) guarantees that the parameter estimate $\hat{\theta}(t)$ is always bounded (even when disturbances are present) and exact estimation in finite time of the unknown parameter $\theta \in \mathbb{R}$ when there exists sufficient excitation of the state (i.e., when there exists $t \geq r$ with $\|x_t\|_2 \geq \sigma$) and the disturbance is absent. Moreover, estimate (2.4) guarantees an ultimate estimation error which is proportional to the magnitude of the disturbance when there exists sufficient excitation of the state: this feature cannot be provided by other adaptive schemes in the literature. The gain of the disturbance for the estimation error (equal to $\sqrt{r}/\sigma$) indicates a trade-off between "learning" and "robustness": when $\sigma > 0$ is large then the gain of the disturbance is small (indicating strong robustness) but this implies that it is less likely the "learning" process to be complete (since it is less likely that there exists $t \geq r$ with $\|x_t\|_2 \geq \sigma$). On the other hand, when $\sigma > 0$ is small then the gain of the disturbance is large (indicating lack of robustness) but this implies that it is more possible the "learning" process to be complete.

## 3. Proof of Theorem 1

We start by stating and proving an auxiliary technical lemma.

**Lemma 1:** *For all* $x \in W^{1,\infty}((-r,0))$ *it holds that*

$$\|x\|_2^2 + \left(\varepsilon^2 - x^2(0)\right)^+ \geq \frac{r\varepsilon^3}{2(3+r)\left(r\sqrt{2}\,\|\dot{x}\|_\infty + \varepsilon\right)} \tag{3.1}$$

*Moreover, for all* $x \in H^1((-r,0))$ *it holds that*

$$\|x\|_2^2 + \left(\varepsilon^2 - x^2(0)\right)^+ \geq \frac{\varepsilon^4 r}{2\varepsilon^2 r + 12\left(\varepsilon^2 + 2r\|\dot{x}\|_2^2\right)} \tag{3.2}$$

**Proof:** Let $x \in W^{1,\infty}((-r,0))$ or $x \in H^1((-r,0))$ be given. We distinguish the following cases.

1st Case: $x^2(0) < \frac{\varepsilon^2}{2}$. In this case we have $\left(\varepsilon^2 - x^2(0)\right)^+ \geq \frac{\varepsilon^2}{2}$ and consequently $\|x\|_2^2 + \left(\varepsilon^2 - x^2(0)\right)^+ \geq \frac{\varepsilon^2}{2}$.

2nd Case: $x^2(0) \geq \frac{\varepsilon^2}{2}$. Using the triangle inequality we get for all $s \in [-r,0]$:

$$\begin{aligned} &|x(0)| \geq \frac{\varepsilon}{\sqrt{2}} \Rightarrow |x(0)| - |x(0)-x(s)| \geq \frac{\varepsilon}{\sqrt{2}} - |x(0)-x(s)| \\ &\Rightarrow |x(s)| \geq \left(\frac{\varepsilon}{\sqrt{2}} - \|\dot{x}\|_\infty |s|\right)^+ \end{aligned} \tag{3.3}$$

and

$$
\begin{aligned}
&|x(0)| \geq \frac{\varepsilon}{\sqrt{2}} \Rightarrow |x(0)| - |x(0)-x(s)| \geq \frac{\varepsilon}{\sqrt{2}} - |x(0)-x(s)| \\
&\Rightarrow |x(s)| \geq \left( \frac{\varepsilon}{\sqrt{2}} - \|\dot{x}\|_2 \sqrt{|s|} \right)^+
\end{aligned}
\tag{3.4}
$$

If $\frac{\varepsilon}{r\sqrt{2}} \geq \|\dot{x}\|_\infty$ then it follows that $\left( \frac{\varepsilon}{\sqrt{2}} - \|\dot{x}\|_\infty |s| \right)^+ = \frac{\varepsilon}{\sqrt{2}} + s\|\dot{x}\|_\infty \geq 0$ for all $s \in [-r,0]$. Thus, in this case we obtain from (3.3):

$$
|x(s)|^2 \geq \frac{\varepsilon^2}{2} + \|\dot{x}\|_\infty^2 s^2 + \sqrt{2}\varepsilon \|\dot{x}\|_\infty s, \forall s \in [-r,0] \Rightarrow \|x\|_2^2 \geq \frac{\varepsilon^2 r}{2} + \frac{1}{3}\|\dot{x}\|_\infty^2 r^3 - \frac{1}{\sqrt{2}}\varepsilon \|\dot{x}\|_\infty r^2
$$

Since the minimum of the function $f(t) = \frac{r^3}{3}t^2 - \frac{\varepsilon r^2}{\sqrt{2}}t + \frac{\varepsilon^2 r}{2}$ for $t \in \left[0, \frac{\varepsilon}{r\sqrt{2}}\right]$ is equal to $\frac{\varepsilon^2 r}{6}$, the above estimate shows that $\|x\|_2^2 \geq \frac{\varepsilon^2 r}{6}$.

If $\frac{\varepsilon}{\sqrt{2r}} \geq \|\dot{x}\|_2$ then it follows that $\left( \frac{\varepsilon}{\sqrt{2}} - \|\dot{x}\|_2 \sqrt{|s|} \right)^+ = \frac{\varepsilon}{\sqrt{2}} - \|\dot{x}\|_2 \sqrt{-s} \geq 0$ for all $s \in [-r,0]$. Thus, in this case we obtain from (3.4):

$$
\begin{aligned}
&|x(s)|^2 \geq \frac{\varepsilon^2}{2} - \|\dot{x}\|_2^2 s - \sqrt{2}\varepsilon \|\dot{x}\|_2 \sqrt{-s}, \forall s \in [-r,0] \\
&\Rightarrow \|x\|_2^2 \geq \frac{\varepsilon^2 r}{2} + \frac{1}{2}\|\dot{x}\|_2^2 r^2 - \frac{2\sqrt{2}}{3}\varepsilon \|\dot{x}\|_2 r\sqrt{r}
\end{aligned}
$$

Since the minimum of the function $f(t) = \frac{r^2}{2}t^2 - \frac{2\sqrt{2r}\varepsilon r}{3}t + \frac{\varepsilon^2 r}{2}$ for $t \in \left[0, \frac{\varepsilon}{\sqrt{2r}}\right]$ is equal to $\frac{\varepsilon^2 r}{12}$, the above estimate shows that $\|x\|_2^2 \geq \frac{\varepsilon^2 r}{12}$.

If $\frac{\varepsilon}{r\sqrt{2}} < \|\dot{x}\|_\infty$ then it follows that $\left( \frac{\varepsilon}{\sqrt{2}} - \|\dot{x}\|_\infty |s| \right)^+ = \frac{\varepsilon}{\sqrt{2}} + s\|\dot{x}\|_\infty \geq 0$ for all $s \in \left[ -\frac{\varepsilon}{\|\dot{x}\|_\infty \sqrt{2}}, 0 \right]$ and $\left( \frac{\varepsilon}{\sqrt{2}} - \|\dot{x}\|_\infty |s| \right)^+ = 0$ for all $s \in \left[ -r, -\frac{\varepsilon}{\|\dot{x}\|_\infty \sqrt{2}} \right]$. Thus, in this case we obtain from (3.3):

$$\left.\begin{array}{c}|x(s)|^2 \geq \frac{\varepsilon^2}{2}+\|\dot{x}\|_\infty^2 s^2+\sqrt{2}\varepsilon\|\dot{x}\|_\infty s, \forall s \in\left[-\frac{\varepsilon}{\|\dot{x}\|_\infty \sqrt{2}}, 0\right] \\ |x(s)|^2 \geq 0, \forall s \in\left[-r,-\frac{\varepsilon}{\|\dot{x}\|_\infty \sqrt{2}}\right]\end{array}\right\} \Rightarrow\|x\|_2^2 \geq \frac{\varepsilon^3}{6 \sqrt{2}\|\dot{x}\|_\infty}$$

If $\frac{\varepsilon}{\sqrt{2r}}<\|\dot{x}\|_2$ then it follows that $\left(\frac{\varepsilon}{\sqrt{2}}-\|\dot{x}\|_2 \sqrt{|s|}\right)^{+}=\frac{\varepsilon}{\sqrt{2}}-\|\dot{x}\|_2 \sqrt{|s|} \geq 0$ for all $s \in\left[-\frac{\varepsilon^2}{2\|\dot{x}\|_2^2}, 0\right]$ and $\left(\frac{\varepsilon}{\sqrt{2}}-\|\dot{x}\|_2 \sqrt{|s|}\right)^{+}=0$ for all $s \in\left[-r,-\frac{\varepsilon^2}{2\|\dot{x}\|_2^2}\right]$. Thus, in this case we obtain from (3.4):

$$\left.\begin{array}{c}|x(s)|^2 \geq \frac{\varepsilon^2}{2}+\|\dot{x}\|_2^2|s|-\sqrt{2}\varepsilon\|\dot{x}\|_2 \sqrt{|s|}, \forall s \in\left[-\frac{\varepsilon^2}{2\|\dot{x}\|_2^2}, 0\right] \\ |x(s)|^2 \geq 0, \forall s \in\left[-r,-\frac{\varepsilon^2}{2\|\dot{x}\|_2^2}\right]\end{array}\right\} \Rightarrow\|x\|_2^2 \geq \frac{\varepsilon^4}{24\|\dot{x}\|_2^2}$$

Combining both cases, we get $\|x\|_2^2 \geq \frac{\varepsilon^2}{6} \min \left(r, \frac{\varepsilon}{\sqrt{2}\|\dot{x}\|_\infty}\right)$, $\|x\|_2^2 \geq \frac{\varepsilon^2}{12} \min \left(\frac{\varepsilon^2}{2\|\dot{x}\|_2^2}, r\right)$ when $\|\dot{x}\|_\infty>0$ (which also implies $\|\dot{x}\|_2>0$) and $\|x\|_2^2 \geq \frac{\varepsilon^2 r}{6}$ when $\|\dot{x}\|_\infty=0$. Thus, we get $\|x\|_2^2+\left(\varepsilon^2-x^2(0)\right)^{+} \geq \frac{\varepsilon^2}{6} \min \left(r, \frac{\varepsilon}{\sqrt{2}\|\dot{x}\|_\infty}\right)$, $\|x\|_2^2+\left(\varepsilon^2-x^2(0)\right)^{+} \geq \frac{\varepsilon^2}{12} \min \left(\frac{\varepsilon^2}{2\|\dot{x}\|_2^2}, r\right)$ when $\|\dot{x}\|_\infty>0$ and $\|x\|_2^2+\left(\varepsilon^2-x^2(0)\right)^{+} \geq \frac{\varepsilon^2 r}{6}$ when $\|\dot{x}\|_\infty=0$. Using the fact that $\min(a, b) \geq \frac{ab}{a+b}$ that holds for all $a, b>0$, we get that $\|x\|_2^2+\left(\varepsilon^2-x^2(0)\right)^{+} \geq \frac{\varepsilon^3 r}{6\left(r \sqrt{2}\|\dot{x}\|_\infty+\varepsilon\right)}$, $\|x\|_2^2+\left(\varepsilon^2-x^2(0)\right)^{+} \geq \frac{\varepsilon^4 r}{12\left(\varepsilon^2+2 r\|\dot{x}\|_2^2\right)}$ when $\|\dot{x}\|_\infty>0$. The inequalities $\|x\|_2^2+\left(\varepsilon^2-x^2(0)\right)^{+} \geq \frac{\varepsilon^3 r}{6\left(r \sqrt{2}\|\dot{x}\|_\infty+\varepsilon\right)}$, $\|x\|_2^2+\left(\varepsilon^2-x^2(0)\right)^{+} \geq \frac{\varepsilon^4 r}{12\left(\varepsilon^2+2 r\|\dot{x}\|_2^2\right)}$ hold even in the case $\|\dot{x}\|_\infty=0$. Consequently, we conclude that the inequalities $\|x\|_2^2+\left(\varepsilon^2-x^2(0)\right)^{+} \geq \frac{\varepsilon^3 r}{6\left(r \sqrt{2}\|\dot{x}\|_\infty+\varepsilon\right)}$, $\|x\|_2^2+\left(\varepsilon^2-x^2(0)\right)^{+} \geq \frac{\varepsilon^4 r}{12\left(\varepsilon^2+2 r\|\dot{x}\|_2^2\right)}$ hold.

Combining the estimate obtained for the 1st Case and the estimates obtained for the 2nd Case, we get:

$$\|x\|_2^2 + \left(\varepsilon^2 - x^2(0)\right)^+ \geq \frac{\varepsilon^2}{2}\min\left(\frac{\varepsilon r}{3\left(r\sqrt{2}\,\|\dot{x}\|_\infty + \varepsilon\right)}, 1\right)$$

$$\|x\|_2^2 + \left(\varepsilon^2 - x^2(0)\right)^+ \geq \frac{\varepsilon^2}{2}\min\left(\frac{\varepsilon^2 r}{6\left(\varepsilon^2 + 2r\|\dot{x}\|_2^2\right)}, 1\right)$$

Using the fact that $\min(a,b) \geq \dfrac{ab}{a+b}$ that holds for all $a,b>0$, we get that

$$\|x\|_2^2 + \left(\varepsilon^2 - x^2(0)\right)^+ \geq \frac{1}{2}\frac{\varepsilon^3 r}{\varepsilon(r+3) + 3r\sqrt{2}\,\|\dot{x}\|_\infty}$$

$$\|x\|_2^2 + \left(\varepsilon^2 - x^2(0)\right)^+ \geq \frac{\varepsilon^4 r}{2\varepsilon^2 r + 12\left(\varepsilon^2 + 2r\|\dot{x}\|_2^2\right)}$$

The above inequalities directly imply inequalities (3.1), (3.2). The proof is complete. $\triangleleft$

We next show the following technical result.

**Lemma 2:** *For every $R \geq 0$ there exists $L(R) > 0$ such that for every $(x,u) \in H^1\left((-r,0)\right) \times L^1\left((-r,0)\right)$, $(y,w) \in H^1\left((-r,0)\right) \times L^1\left((-r,0)\right)$ with $\max\left(\|x\|_\infty, \|\dot{x}\|_2, \|u\|_1\right) \leq R$, $\max\left(\|y\|_\infty, \|\dot{y}\|_2, \|w\|_1\right) \leq R$, it holds that*

$$\left|p(x,u) - p(y,w)\right| \leq L(R)\left(\|x-y\|_\infty + \|u-w\|_1\right) \tag{3.5}$$

**Proof:** Let $R \geq 0$, $(x,u) \in H^1\left((-r,0)\right) \times L^1\left((-r,0)\right)$, $(y,w) \in H^1\left((-r,0)\right) \times L^1\left((-r,0)\right)$ with $\max\left(\|x\|_\infty, \|\dot{x}\|_2, \|u\|_1\right) \leq R$, $\max\left(\|y\|_\infty, \|\dot{y}\|_2, \|w\|_1\right) \leq R$ be given. Using definition (2.7) and the triangle inequality we get:

$$\begin{aligned}
\left|p(x,u) - p(y,w)\right| &\leq \frac{cr}{2\left(\|x\|_2^2 + \left(\varepsilon - x^2(0)\right)^+\right)}\left|x^2(0) - y^2(0)\right| \\
&+ \frac{1}{2\left(\|x\|_2^2 + \left(\varepsilon - x^2(0)\right)^+\right)}\left|\left(x^2(0) - x^2(-r) - 2\langle x,u\rangle\right)^+ - \left(y^2(0) - y^2(-r) - 2\langle y,w\rangle\right)^+\right| \\
&+ \frac{\left(y^2(0) - y^2(-r) - 2\langle y,w\rangle\right)^+ + cry^2(0)}{2\left(\|x\|_2^2 + \left(\varepsilon - x^2(0)\right)^+\right)\left(\|y\|_2^2 + \left(\varepsilon - y^2(0)\right)^+\right)}\left|\|y\|_2^2 - \|x\|_2^2\right| \\
&+ \frac{\left(y^2(0) - y^2(-r) - 2\langle y,w\rangle\right)^+ + cry^2(0)}{2\left(\|x\|_2^2 + \left(\varepsilon - x^2(0)\right)^+\right)\left(\|y\|_2^2 + \left(\varepsilon - y^2(0)\right)^+\right)}\left|\left(\varepsilon - y^2(0)\right)^+ - \left(\varepsilon - x^2(0)\right)^+\right|
\end{aligned} \tag{3.6}$$

Using (3.2), the fact that $\left|(a)^{+}-(b)^{+}\right|\leq |a-b|$ for all $a,b\in\mathbb{R}$, the fact that $\max\left(\|\dot{x}\|_2,\|\dot{y}\|_2\right)\leq R$ and the triangle inequality, we obtain from (3.6):

$$\begin{aligned}&|p(x,u)-p(y,w)|\leq (1+cr)H(R)\left|x^2(0)-y^2(0)\right|\\&+H(R)\left|x^2(-r)-y^2(-r)\right|+2H(R)\left|\langle x,u\rangle-\langle y,w\rangle\right|\\&+2H^2(R)\left(\left(y^2(0)-y^2(-r)-2\langle y,w\rangle\right)^{+}+cry^2(0)\right)\left|\|y\|_2^2-\|x\|_2^2\right|\\&+2H^2(R)\left(\left(y^2(0)-y^2(-r)-2\langle y,w\rangle\right)^{+}+cry^2(0)\right)\left|x^2(0)-y^2(0)\right|\end{aligned}\tag{3.7}$$

where

$$H(R):=\frac{\varepsilon^2 r+6\left(\varepsilon^2+2rR^2\right)}{\varepsilon^4 r}\tag{3.8}$$

Using the fact that $\left(y^2(0)-y^2(-r)-2\langle y,w\rangle\right)^{+}\leq y^2(0)+2\left|\langle y,w\rangle\right|$ in conjunction with (3.7) and the triangle inequality, we get:

$$\begin{aligned}&|p(x,u)-p(y,w)|\leq (1+cr)H(R)\left(|x(0)|+|y(0)|\right)|x(0)-y(0)|\\&+H(R)\left(|x(-r)|+|y(-r)|\right)|x(-r)-y(-r)|+2H(R)\left|\langle x,u-w\rangle\right|\\&+2H(R)\left|\langle x-y,w\rangle\right|\\&+2H^2(R)\left((1+cr)y^2(0)+2\left|\langle y,w\rangle\right|\right)\left(\|y\|_2+\|x\|_2\right)\left|\|y\|_2-\|x\|_2\right|\\&+2H^2(R)\left((1+cr)y^2(0)+2\left|\langle y,w\rangle\right|\right)\left(|x(0)|+|y(0)|\right)|x(0)-y(0)|\end{aligned}\tag{3.9}$$

Using (3.9), the fact that $\max\left(\|x\|_\infty,\|y\|_\infty\right)\leq R$, Hölder's inequality and the triangle inequality, we get:

$$\begin{aligned}&|p(x,u)-p(y,w)|\leq 2R(1+cr)H(R)|x(0)-y(0)|\\&+2RH(R)|x(-r)-y(-r)|+2H(R)\|x\|_\infty\|u-w\|_1\\&+2H(R)\|x-y\|_\infty\|w\|_1\\&+2H^2(R)\left((1+cr)R^2+2\|y\|_\infty\|w\|_1\right)\left(\|y\|_2+\|x\|_2\right)\|x-y\|_2\\&+4RH^2(R)\left((1+cr)R^2+2\|y\|_\infty\|w\|_1\right)|x(0)-y(0)|\end{aligned}\tag{3.10}$$

Using the fact that $\|f\|_2\leq\sqrt{r}\,\|f\|_\infty$ that holds for all $f\in C^0\left([-r,0];\mathbb{R}\right)$, we get from (3.10):

$$\begin{aligned}&|p(x,u)-p(y,w)|\leq 2R(1+cr)H(R)\|x-y\|_\infty\\&+2RH(R)\|x-y\|_\infty+2H(R)\|x\|_\infty\|u-w\|_1\\&+2H(R)\|x-y\|_\infty\|w\|_1\\&+2rH^2(R)\left((1+cr)R^2+2\|y\|_\infty\|w\|_1\right)\left(\|y\|_\infty+\|x\|_\infty\right)\|x-y\|_\infty\\&+4RH^2(R)\left((1+cr)R^2+2\|y\|_\infty\|w\|_1\right)\|x-y\|_\infty\end{aligned}\tag{3.11}$$

Using the fact that $\max\left(\|x\|_\infty, \|y\|_\infty, \|w\|_1\right) \le R$ and (3.11) we obtain (3.5) with

$$L(R) = (3+cr)\left(1+2R^2(r+1)H(R)\right)2RH(R)$$

The proof is complete. $\triangleleft$

Now we are ready to provide the first existence result for system (2.1), (2.3).

**Lemma 3:** *Given* $(x_0,u_0)\in X$ *and* $d\in L^\infty(\mathbb{R}_+)$ *there exist* $\delta>0$ *and continuous mappings* $x,u:[-r,\delta]\to\mathbb{R}$ *with* $x\in W^{1,\infty}((-r,\delta))$, $x(s)=x_0(s)$, $u(s)=u_0(s)$ *for all* $s\in[-r,0]$ *such that (2.1) holds for* $t\in[0,\delta]$ *a.e. and (2.3) holds for all* $t\in[0,\delta]$.

**Proof:** Let arbitrary $(x_0,u_0)\in X$ and $d\in L^\infty(\mathbb{R}_+)$ be given. Define the closed set

$$A:=\left\{ a\in C^0([0,\delta];\mathbb{R}_+): a(0)=p(x_0,u_0), \max_{0\le s\le\delta}(a(s))\le M \right\}\subset C^0([0,\delta]) \tag{3.12}$$

where $\delta, M>0$ are constants (to be selected appropriately). We also define the mapping $P:A\to C^0([0,\delta];\mathbb{R})$ in the following way. Given $a\in A$ we define

$$y(t)=\exp\left((\theta-2c)t-\int_0^t a(s)ds\right)x_0(0)+\int_0^t \exp\left((\theta-2c)(t-\tau)-\int_\tau^t a(s)ds\right)d(\tau)d\tau,$$
$$\text{for } t\in(0,\delta] \tag{3.13}$$

$$y(s)=x_0(s), \text{ for all } s\in[-r,0] \tag{3.14}$$

$$w(t)=-(2c+a(t))y(t), \text{ for } t\in(0,\delta] \tag{3.15}$$

$$w(s)=u_0(s), \text{ for all } s\in[-r,0] \tag{3.16}$$

$$\bar{a}(t)=p(y_t,w_t)=\frac{\left(y^2(t)-y^2(t-r)-2\langle y_t,w_t\rangle\right)^+ + cry^2(t)}{2\left(\|y_t\|_2^2+\left(\varepsilon-y^2(t)\right)^+\right)}, \text{ for } t\in[0,\delta] \tag{3.17}$$

We define:

$$P(a)=\bar{a} \tag{3.18}$$

In order to complete the proof of the theorem, it suffices to show that the mapping $P:A\to C^0([0,\delta];\mathbb{R})$ has a fixed point. Therefore, by using Banach's fixed point theorem, we have to show that there exist $\delta, M>0$ such that the mapping takes values in $A$ and is a contraction.

It should be noticed that definitions (3.12), (3.14), (3.16) and (3.17) guarantee that $\bar{a}(0)=p(x_0,u_0)$. Moreover, definitions (3.12), (3.13), (3.14), (3.15), (3.16) guarantee that

$$\max_{-r\le t\le\delta}\left(|y(t)|\right)\le\exp\left((\theta-2c)^+\delta\right)\left(\|x_0\|_\infty+\delta\|d\|_\infty\right) \tag{3.19}$$

$$\dot{y}(t)=(\theta-2c-a(t))y(t)+d(t)\text{, for } t\in(0,\delta] \text{ a.e.} \tag{3.20}$$

$$\sup_{0<t<\delta}\left(|\dot{y}(t)|\right)\le\left(|\theta-2c|+M\right)\exp\left((\theta-2c)^+\delta\right)\left(\|x_0\|_\infty+\delta\|d\|_\infty\right)+\|d\|_\infty \tag{3.21}$$

$$\begin{aligned}&\int_{-r}^{\delta}\dot{y}^2(s)ds\le\|\dot{x}_0\|_2^2\\&+\delta\left(\|d\|_\infty+\left(|\theta-2c|+M\right)\exp\left((\theta-2c)^+\delta\right)\left(\|x_0\|_\infty+\delta\|d\|_\infty\right)\right)^2\end{aligned} \tag{3.22}$$

$$\max_{0\le t\le\delta}\left(|w(t)|\right)\le\left(2c+M\right)\exp\left((\theta-2c)^+\delta\right)\left(\|x_0\|_\infty+\delta\|d\|_\infty\right) \tag{3.23}$$

$$\int_{-r}^{\delta}|w(s)|ds\le\|u_0\|_1+\delta\left(2c+M\right)\exp\left((\theta-2c)^+\delta\right)\left(\|x_0\|_\infty+\delta\|d\|_\infty\right) \tag{3.24}$$

Using definition (3.17), inequality (3.2), the fact that $\left(y^2(t)-y^2(t-r)-2\langle y_t,w_t\rangle\right)^+\le y^2(t)+2\left|\langle y_t,w_t\rangle\right|$, Hölder's inequality (which implies that $\left|\langle y_t,w_t\rangle\right|\le\|y_t\|_\infty\|w_t\|_1$) and estimates (3.19), (3.22), (3.24), we obtain the estimate for all $\delta\in[0,1]$:

$$\bar{a}(t)\le gH+\left(\lambda H+2R(g+\lambda)\left(2c+M\right)\right)\delta\text{, for all } t\in[0,\delta] \tag{3.25}$$

where

$$\begin{aligned}R&=\exp\left((\theta-2c)^+\right)\left(\|x_0\|_\infty+\|d\|_\infty\right)\\g&=\frac{12R}{\varepsilon^4}\left(\frac{\varepsilon^2(r+6)}{12r}+\|\dot{x}_0\|_2^2\right)\\H&=\left((1+cr)R+2\|u_0\|_1\right)\\\lambda&=\frac{12R}{\varepsilon^4}\left(\|d\|_\infty+|\theta-2c|R+MR\right)^2\end{aligned}$$

It follows from (3.12) and (3.25) that for $M=gH+1$ there exists $\delta>0$ sufficiently small so that the mapping $P:A\to C^0\left([0,\delta];\mathbb{R}\right)$ takes values in $A$.

Given $\tilde{a}\in A$ we define

$$z(t)=\exp\left((\theta-2c)t-\int_0^t\tilde{a}(s)ds\right)x_0(0)+\int_0^t\exp\left((\theta-2c)(t-\tau)-\int_\tau^t\tilde{a}(s)ds\right)d(\tau)d\tau\text{,}$$
$$\text{for } t\in(0,\delta] \tag{3.26}$$

$$z(s)=x_0(s)\text{, for all } s\in[-r,0] \tag{3.27}$$

$$v(t) = -\left(2c + \tilde{a}(t)\right) z(t) \text{, for } t \in (0,\delta] \tag{3.28}$$

$$v(s) = u_0(s) \text{, for all } s \in [-r,0] \tag{3.29}$$

$$\left(P(\tilde{a})\right)(t) = p(z_t, v_t) \text{, for } t \in [0,\delta] \tag{3.30}$$

Using (3.13) and (3.26), the fact that $\left|\exp(-\varphi) - \exp(-\eta)\right| \le \left|\varphi - \eta\right|$ for all $\varphi, \eta \ge 0$, we get for all $t \in [0,\delta]$:

$$\begin{aligned}
\left|y(t) - z(t)\right| &\le \exp\left((\theta - 2c)t\right)\left|\exp\left(-\int_0^t a(s)ds\right) - \exp\left(-\int_0^t \tilde{a}(s)ds\right)\right|\left\|x_0\right\|_\infty \\
&+ \left\|d\right\|_\infty \int_0^t \exp\left((\theta - 2c)(t-\tau)\right)\left|\exp\left(-\int_\tau^t a(s)ds\right) - \exp\left(-\int_\tau^t \tilde{a}(s)ds\right)\right| d\tau \\
&\le \exp\left((\theta - 2c)^+ \delta\right)\left|\int_0^t \left(a(s) - \tilde{a}(s)\right)ds\right|\left\|x_0\right\|_\infty \\
&+ \left\|d\right\|_\infty \exp\left((\theta - 2c)^+ \delta\right)\int_0^t \left|\int_\tau^t \left(a(s) - \tilde{a}(s)\right)ds\right| d\tau \\
&\le \exp\left((\theta - 2c)^+ \delta\right)\left(\left\|x_0\right\|_\infty + \left\|d\right\|_\infty \delta\right)\delta \max_{0 \le s \le \delta}\left(\left|a(s) - \tilde{a}(s)\right|\right)
\end{aligned}$$

Combining (3.14), (3.27) and the above estimate we get for all $t \in [0,\delta]$:

$$\left\|y_t - z_t\right\|_\infty \le \exp\left((\theta - 2c)^+ \delta\right)\left(\left\|x_0\right\|_\infty + \left\|d\right\|_\infty \delta\right)\delta \max_{0 \le s \le \delta}\left(\left|a(s) - \tilde{a}(s)\right|\right) \tag{3.31}$$

Using (3.12), (3.15), (3.19), (3.28) and (3.31) we obtain for all $t \in [0,\delta]$:

$$\left|w(t) - v(t)\right| \le \exp\left((\theta - 2c)^+ \delta\right)\left((2c + M)\delta + 1\right)\left(\left\|x_0\right\|_\infty + \left\|d\right\|_\infty \delta\right)\max_{0 \le s \le \delta}\left(\left|a(s) - \tilde{a}(s)\right|\right) \tag{3.32}$$

Consequently, we get from (3.16), (3.29) and (3.32) for all $t \in [0,\delta]$:

$$\left\|w_t - v_t\right\|_1 \le \delta \exp\left((\theta - 2c)^+ \delta\right)\left((2c + M)\delta + 1\right)\left(\left\|x_0\right\|_\infty + \left\|d\right\|_\infty \delta\right)\max_{0 \le s \le \delta}\left(\left|a(s) - \tilde{a}(s)\right|\right) \tag{3.33}$$

It should be noticed that definitions (3.26)-(3.29) guarantee that all estimates (3.19), (3.21), (3.22), (3.23), (3.24) hold with $z, v$ in place of $y, w$. Therefore, for $\delta \in [0,1]$ it follows that $\max\left(\left\|z_t\right\|_\infty, \left\|\dot{z}_t\right\|_2, \left\|v_t\right\|_1\right) \le \bar{R}$, $\max\left(\left\|y_t\right\|_\infty, \left\|\dot{y}_t\right\|_2, \left\|w_t\right\|_1\right) \le \bar{R}$ for all $t \in [0,\delta]$ with

$$\begin{aligned}
\bar{R} &= \exp\left((\theta - 2c)^+\right)\left(\left\|x_0\right\|_\infty + \left\|d\right\|_\infty\right) + \left\|\dot{x}_0\right\|_2^2 \\
&+ \left(\left\|d\right\|_\infty + \left(\left|\theta - 2c\right| + M\right)\exp\left((\theta - 2c)^+\right)\left(\left\|x_0\right\|_\infty + \left\|d\right\|_\infty\right)\right)^2 \\
&+ \left\|u_0\right\|_1 + \left(2c + M\right)\exp\left((\theta - 2c)^+\right)\left(\left\|x_0\right\|_\infty + \left\|d\right\|_\infty\right)
\end{aligned} \tag{3.34}$$

Lemma 2 implies that there exists $L(\bar{R})>0$ such that for every $(x,u)\in H^1\left((-r,0)\right)\times L^1\left((-r,0)\right)$, $(y,w)\in H^1\left((-r,0)\right)\times L^1\left((-r,0)\right)$ with $\max\left(\|x\|_\infty,\|\dot{x}\|_2,\|u\|_1\right)\le \bar{R}$, $\max\left(\|y\|_\infty,\|\dot{y}\|_2,\|w\|_1\right)\le \bar{R}$, it holds that

$$\left|p(x,u)-p(y,w)\right|\le L(\bar{R})\left(\|x-y\|_\infty+\|u-w\|_1\right) \tag{3.35}$$

Definitions (3.17), (3.18), (3.30) and inequality (3.35) imply that for all $\delta\in[0,1]$ it holds that

$$\left|(P(a))(t)-(P(\tilde{a}))(t)\right|\le L(\bar{R})\left(\|y_t-z_t\|_\infty+\|w_t-v_t\|_1\right) \tag{3.36}$$

Combining (3.31), (3.33) and (3.36), we obtain for all $\delta\in[0,1]$:

$$\begin{aligned}&\max_{0\le t\le\delta}\left(\left|(P(a))(t)-(P(\tilde{a}))(t)\right|\right)\\ &\le \delta L(\bar{R})\exp\left((\theta-2c)^+\right)\left(\|x_0\|_\infty+\|d\|_\infty\right)\left(2c+M+2\right)\max_{0\le s\le\delta}\left(\left|a(s)-\tilde{a}(s)\right|\right)\end{aligned} \tag{3.37}$$

It follows from (3.37) that there exists $\delta>0$ sufficiently small so that the mapping $P:A\to C^0\left([0,\delta];\mathbb{R}\right)$ is a contraction. The proof is complete. $\lhd$

The following result deals with continuity with respect to initial conditions for the solutions for system (2.1), (2.3). It also shows uniqueness of solutions for system (2.1), (2.3).

**Lemma 4:** *For every $T,R>0$, there exists $Q(R,T)>0$ with the following property:*

**(P)** *For every $d\in L^\infty(\mathbb{R}_+)$ with $\|d\|_\infty\le R$, $x,u,y,w\in C^0\left([-r,T]\right)$ with $x,y\in W^{1,\infty}\left((-r,T)\right)$ and $\max_{-r\le s\le T}\left(|x(s)|\right)\le R$, $\max_{-r\le s\le T}\left(|y(s)|\right)\le R$, $\max_{-r\le s\le T}\left(|u(s)|\right)\le R$, $\max_{-r\le s\le T}\left(|w(s)|\right)\le R$, $\sup_{-r<s<T}\left(|\dot{x}(s)|\right)\le R$, $\sup_{-r<s<T}\left(|\dot{y}(s)|\right)\le R$, for which (2.1) holds for $t\in[0,T]$ a.e., (2.3) holds for all $t\in[0,T]$, $\dot{y}(t)=\theta y(t)+w(t)+d(t)$ holds for $t\in[0,T]$ a.e., $w(t)=-\left(2c+p\left(y_t,w_t\right)\right)y(t)$ holds for all $t\in[0,T]$, the following estimate holds:*

$$\|x_t-y_t\|_\infty+\|u_t-w_t\|_1\le Q(R,T)\left(\|x_0-y_0\|_\infty+\|u_0-w_0\|_1\right),\text{ for all } t\in[0,T] \tag{3.38}$$

**Proof:** Let arbitrary $T,R>0$ be given. By virtue of Lemma 2 there exists $L(R)>0$ such that for every $(x,u)\in H^1\left((-r,0)\right)\times L^1\left((-r,0)\right)$, $(y,w)\in H^1\left((-r,0)\right)\times L^1\left((-r,0)\right)$ with $\max\left(\|x\|_\infty,\|\dot{x}\|_2,\|u\|_1\right)\le R$, $\max\left(\|y\|_\infty,\|\dot{y}\|_2,\|w\|_1\right)\le R$, (3.5) holds.

Let arbitrary $d\in L^\infty(\mathbb{R}_+)$ with $\|d\|_\infty\le R$, $x,u,y,w\in C^0\left([-r,T]\right)$ with $x,y\in W^{1,\infty}\left((-r,T)\right)$ and $\max_{-r\le s\le T}\left(|x(s)|\right)\le R$, $\max_{-r\le s\le T}\left(|y(s)|\right)\le R$, $\max_{-r\le s\le T}\left(|u(s)|\right)\le R$, $\max_{-r\le s\le T}\left(|w(s)|\right)\le R$, $\sup_{-r<s<T}\left(|\dot{x}(s)|\right)\le R$, $\sup_{-r<s<T}\left(|\dot{y}(s)|\right)\le R$, for which (2.1) holds for $t\in[0,T]$ a.e., (2.3) holds for all $t\in[0,T]$,

$\dot{y}(t)=\theta y(t)+w(t)+d(t)$ holds for $t\in[0,T]$ a.e., $w(t)=-\left(2c+p\left(y_t,w_t\right)\right)y(t)$ holds for all $t\in[0,T]$, be given. Since $\|\varphi\|_2\leq\sqrt{r}\,\|\varphi\|_\infty$, $\|\varphi\|_1\leq r\,\|\varphi\|_\infty$ for all $\varphi\in C^0\left([-r,0]\right)$ we get from (2.3), (3.5) and the equation $w(t)=-\left(2c+p\left(y_t,w_t\right)\right)y(t)$ for all $t\in[0,T]$:

$$
\begin{aligned}
&|u(t)-w(t)|\leq\left(2c+p\left(x_t,u_t\right)\right)|x(t)-y(t)|+|y(t)|\left|p\left(x_t,u_t\right)-p\left(y_t,w_t\right)\right|\\
&\leq\left(2c+p\left(x_t,u_t\right)\right)|x(t)-y(t)|\\
&+|y(t)|\,L\left((1+r)R\right)\|x_t-y_t\|_\infty+|y(t)|\,L\left((1+r)R\right)\|u_t-w_t\|_1\\
&\leq\left(2c+p\left(x_t,u_t\right)-p\left(0,0\right)\right)|x(t)-y(t)|\\
&+RL\left((1+r)R\right)\|x_t-y_t\|_\infty+RL\left((1+r)R\right)\|u_t-w_t\|_1\\
&\leq\left(2c+L\left((1+r)R\right)\|x_t\|_\infty+L\left((1+r)R\right)\|u_t\|_1\right)|x(t)-y(t)|\\
&+RL\left((1+r)R\right)\|x_t-y_t\|_\infty+RL\left((1+r)R\right)\|u_t-w_t\|_1\\
&\leq\left(2c+(2+r)RL\left((1+r)R\right)\right)\left(\|x_t-y_t\|_\infty+\|u_t-w_t\|_1\right)
\end{aligned}
\tag{3.39}
$$

Consequently, we obtain from (3.39) for all $t\in[0,T]$:

$$
\|u_t-w_t\|_1\leq\|u_0-w_0\|_1+\left(2c+(2+r)RL\left((1+r)R\right)\right)\int_0^t\left(\|x_s-y_s\|_\infty+\|u_s-w_s\|_1\right)ds \tag{3.40}
$$

Using (2.1) and the equation $\dot{y}(t)=\theta y(t)+w(t)+d(t)$ we get for all $t\in[0,T]$:

$$
|x(t)-y(t)|\leq|x(0)-y(0)|+|\theta|\int_0^t|x(s)-y(s)|\,ds+\int_0^t|u(s)-w(s)|\,ds \tag{3.41}
$$

Consequently, we obtain from (3.39) and (3.41) for all $t\in[0,T]$:

$$
\|x_t-y_t\|_\infty\leq\|x_0-y_0\|_\infty+\left(2c+|\theta|+(2+r)RL\left((1+r)R\right)\right)\int_0^t\left(\|x_s-y_s\|_\infty+\|u_s-w_s\|_1\right)ds \tag{3.42}
$$

Combining (3.40) and (3.42), we get for all $t\in[0,T]$:

$$
b(t)\leq b(0)+2\left(2c+|\theta|+(2+r)RL\left((1+r)R\right)\right)\int_0^t b(s)ds \tag{3.43}
$$

where $b(t)=\|x_t-y_t\|_\infty+\|u_t-w_t\|_1$. Estimate (3.38) with $Q(R,T)=\exp\left(2\left(2c+|\theta|+(2+r)RL\left((1+r)R\right)\right)T\right)$ is a consequence of (3.43) and the Gronwall-Bellman Lemma. The proof is complete. $\triangleleft$

We are now ready to prove Theorem 1.

**Proof of Theorem 1:** Let arbitrary $(x_0,u_0)\in X$ and $d\in L^\infty(\mathbb{R}_+)$ be given.

We say that property $\Pi(\delta)$ holds for some $\delta>0$ if there exist continuous mappings $x,u:[-r,\delta]\to\mathbb{R}$ with $x\in W^{1,\infty}((-r,\delta))$, $x(s)=x_0(s)$, $u(s)=u_0(s)$ for all $s\in[-r,0]$ such that (2.1) holds for $t\in[0,\delta]$ a.e. and (2.3) holds for all $t\in[0,\delta]$. We define:

$$t_{\max}=\sup\left\{\,\delta>0:\Pi(\delta)\ \textit{holds}\,\right\} \tag{3.44}$$

By virtue of Lemma 3, the set $\left\{\,\delta>0:\Pi(\delta)\ \textit{holds}\,\right\}$ is non-empty. Consequently, definition (3.44) implies that $t_{\max}\in(0,+\infty]$. Moreover, definition (3.44) implies that there exist continuous mappings $x,u:[-r,t_{\max})\to\mathbb{R}$ with $x\in W^{1,\infty}((-r,T))$ for all $T<t_{\max}$, $x(s)=x_0(s)$, $u(s)=u_0(s)$ for all $s\in[-r,0]$ such that (2.1) holds for $t\in[0,t_{\max})$ a.e. and (2.3) holds for all $t\in[0,t_{\max})$. Uniqueness of the pair of mappings $x,u:[-r,t_{\max})\to\mathbb{R}$ follows from Lemma 4.

Next suppose that $t_{\max}<+\infty$. We show next by means of a contradiction argument that $u:[-r,t_{\max})\to\mathbb{R}$ is not bounded. Using (2.1), (2.3) we establish the following estimate:

$$|x(t)|\le\exp\left((\theta-2c)^+t_{\max}\right)\left(\|x_0\|_\infty+t_{\max}\|d\|_\infty\right),\ \text{for all } t\in[0,t_{\max}) \tag{3.45}$$

It follows from (3.45) that $x:[-r,t_{\max})\to\mathbb{R}$ is bounded. If $u:[-r,t_{\max})\to\mathbb{R}$ is bounded, then both mappings $x,u:[-r,t_{\max})\to\mathbb{R}$ can be extended continuously on $[-r,t_{\max}]$. Indeed, notice that since $u:[-r,t_{\max})\to\mathbb{R}$ is bounded, it follows from (2.1) that $\dot{x}(t)$, $\frac{d}{dt}\langle x_t,u_t\rangle$ are also (essentially) bounded on $[-r,t_{\max})$ and $[0,t_{\max})$, respectively. Consequently, the mappings $x:[-r,t_{\max})\to\mathbb{R}$ and $[0,t_{\max})\ni t\mapsto\left(x^2(t)-x^2(t-r)-2\langle x_t,u_t\rangle\right)^+\in\mathbb{R}$ are uniformly continuous. Thus, the mapping $x:[-r,t_{\max})\to\mathbb{R}$ can be extended continuously on $[-r,t_{\max}]$. Lemma 2.1 on page 40 in [5] implies that the mapping $[-r,t_{\max}]\ni t\mapsto x_t\in C^0([-r,0])$ is continuous and consequently, the mapping $[0,t_{\max}]\ni t\mapsto\|x_t\|_2^2+\left(\varepsilon^2-x^2(t)\right)^+\in\mathbb{R}$ is uniformly continuous and positive. Hence definition (2.7) implies that the mapping $[0,t_{\max})\ni t\mapsto p(x_t,u_t)\in\mathbb{R}$ is uniformly continuous and by virtue of (2.3) it follows that the mapping $[-r,t_{\max})\ni t\mapsto u(t)\in\mathbb{R}$ is uniformly continuous. Thus, the mapping $u:[-r,t_{\max})\to\mathbb{R}$ can be extended continuously on $[-r,t_{\max}]$.

Since $u:[-r,t_{\max})\to\mathbb{R}$ is bounded, it follows from (2.1) that $\dot{x}$ is also (essentially) bounded on $[-r,t_{\max})$ and consequently $x_{t_{\max}}\in W^{1,\infty}((-r,0))$. Moreover, by means of (2.1), (2.3), continuity of $p$ (recall Lemma 2), and continuity of the mappings $[-r,t_{\max}]\ni t\mapsto x_t,u_t\in C^0([-r,0])$ (recall Lemma 2.1 on page 40 in [5]) it follows that $u(t_{\max})=-\left(2c+p\left(x_{t_{\max}},u_{t_{\max}}\right)\right)x(t_{\max})$. Consequently, definition (2.9) implies that $\left(x_{t_{\max}},u_{t_{\max}}\right)\in X$. Therefore, by exploiting Lemma 3, we can find $\delta>0$ such that property $\Pi(t_{\max}+\delta)$ holds. This contradicts definition (3.44).

The proof is complete. $\lhd$

## 4. Proofs of Theorem 2 and Theorem 3

**Proof of Theorem 2:** Let $(x_0,u_0)\in X$ and $d\in L^{\infty}(\mathbb{R}_+)$ be given. Then, by virtue of Theorem 1, there exist $t_{\max}\in(0,+\infty]$ and a unique pair of continuous mappings $x,u:[-r,t_{\max})\to\mathbb{R}$ with $x\in W^{1,\infty}((-r,T))$ for all $T<t_{\max}$, $x(s)=x_0(s)$, $u(s)=u_0(s)$ for all $s\in[-r,0]$ such that (2.1) holds for $t\in[0,t_{\max})$ a.e. and (2.3) holds for all $t\in[0,t_{\max})$. Moreover, if $t_{\max}<+\infty$ then $\limsup_{t\to t_{\max}^-}(|u(t)|)=+\infty$.

Defining

$$d(t)=\dot{x}(t)-\theta x(t)-u(t)\text{, for } t\in(-r,0)\text{ a.e.} \tag{4.1}$$

we conclude that (2.1) holds for $t\in[-r,t_{\max})$ a.e..

By virtue of (2.1), we conclude that the following identity is valid for all $t\in[0,t_{\max})$:

$$x^2(t)-x^2(t-r)-2\langle x_t,u_t\rangle=2\theta\|x_t\|_2^2+2\langle x_t,d_t\rangle \tag{4.2}$$

Combining (2.3), (2.7) and (4.2), we get for all $t\in[0,t_{\max})$:

$$u(t)=-2cx(t)-\frac{2\left(\theta\|x_t\|_2^2+\langle x_t,d_t\rangle\right)^+ +crx^2(t)}{2\left(\|x_t\|_2^2+\left(\varepsilon^2-x^2(t)\right)^+\right)}x(t) \tag{4.3}$$

Define the function

$$W(x):=\frac{1}{2}\left(\left(x^2-\varepsilon^2\right)^+\right)^2\text{, for all } x\in\mathbb{R} \tag{4.4}$$

We next assume that $t_{\max}>r$.

Using definition (4.4) and (2.1), (4.3), we get for $t\in[r,t_{\max})$ a.e.:

$$\frac{d}{dt}W(x(t))=-\left(x^2(t)-\varepsilon^2\right)^+2x^2(t)\left(2c+\frac{2\left(\theta\|x_t\|_2^2+\langle x_t,d_t\rangle\right)^+ +crx^2(t)}{2\left(\|x_t\|_2^2+\left(\varepsilon^2-x^2(t)\right)^+\right)}-\theta\right) \tag{4.5}$$

$$+\left(x^2(t)-\varepsilon^2\right)^+2x(t)d(t)$$

Next, we focus on the case $x^2(t)\geq\varepsilon^2$. Using (4.5) and the fact that $2x(t)d(t)\leq cx^2(t)+c^{-1}d^2(t)$ we get for $t\in[r,t_{\max})$ a.e. when $x^2(t)\geq\varepsilon^2$:

$$\frac{d}{dt}W(x(t)) \le -\left(x^2(t)-\varepsilon^2\right)x^2(t)\left(3c+\frac{2\left(\theta\|x_t\|_2^2+\langle x_t,d_t\rangle\right)^+ + crx^2(t)}{\|x_t\|_2^2}-2\theta\right) \tag{4.6}$$
$$+c^{-1}\left(x^2(t)-\varepsilon^2\right)\|d\|_\infty^2$$

Using the fact that $\left(\theta\|x_t\|_2^2+\langle x_t,d_t\rangle\right)^+ \ge (\theta)^+\|x_t\|_2^2-\left(-\langle x_t,d_t\rangle\right)^+$ (a consequence of the fact $(a+b)^+\le(a)^+ + (b)^+$ that holds for all $a,b\in\mathbb{R}$), we get from (4.6) the following inequality for $t\in[r,t_{\max})$ a.e. when $x^2(t)\ge\varepsilon^2$:

$$\frac{d}{dt}W(x(t)) \le -\left(x^2(t)-\varepsilon^2\right)x^2(t)\left(3c+\frac{crx^2(t)-2\left(-\langle x_t,d_t\rangle\right)^+}{\|x_t\|_2^2}+2(\theta)^+-2\theta\right) \tag{4.7}$$
$$+c^{-1}\left(x^2(t)-\varepsilon^2\right)\|d\|_\infty^2$$

The Cauchy-Schwarz inequality implies that $\left|\langle x_t,d_t\rangle\right|\le\|x_t\|_2\|d_t\|_2\le\sqrt{r}\,\|x_t\|_2\|d\|_\infty$. Using the previous estimate in conjunction with (4.7) and the fact that $(\theta)^+\ge\theta$, we get the following inequality for $t\in[r,t_{\max})$ a.e. when $x^2(t)\ge\varepsilon^2$:

$$\frac{d}{dt}W(x(t)) \le -\left(x^2(t)-\varepsilon^2\right)x^2(t)\left(3c+\frac{crx^2(t)}{\|x_t\|_2^2}\right)$$
$$+\left(x^2(t)-\varepsilon^2\right)x^2(t)\frac{2\sqrt{r}\,\|d\|_\infty}{\|x_t\|_2}+c^{-1}\left(x^2(t)-\varepsilon^2\right)\|d\|_\infty^2 \tag{4.8}$$

Using the estimate $x^2(t)\frac{2\sqrt{r}\,\|d\|_\infty}{\|x_t\|_2}\le c^{-1}\|d\|_\infty^2+cr\frac{x^4(t)}{\|x_t\|_2^2}$, we get from (4.8) the following inequality for $t\in[r,t_{\max})$ a.e. when $x^2(t)\ge\varepsilon^2$:

$$\frac{d}{dt}W(x(t)) \le -3c\left(x^2(t)-\varepsilon^2\right)x^2(t)+2c^{-1}\left(x^2(t)-\varepsilon^2\right)\|d\|_\infty^2 \tag{4.9}$$

Re-arranging the terms in (4.9) we obtain the following inequality for $t\in[r,t_{\max})$ a.e. when $x^2(t)\ge\varepsilon^2$:

$$\frac{d}{dt}W(x(t)) \le -3c\left(x^2(t)-\varepsilon^2\right)^2$$
$$-3c\varepsilon^2\left(x^2(t)-\varepsilon^2\right)+2c^{-1}\left(x^2(t)-\varepsilon^2\right)\|d\|_\infty^2 \tag{4.10}$$

Using the inequality $2\left(x^2(t)-\varepsilon^2\right)\|d\|_\infty^2 \le c^2\left(x^2(t)-\varepsilon^2\right)^2 + c^{-2}\|d\|_\infty^4$ we get from (4.10) the following inequality for $t\in[r,t_{\max})$ a.e. when $x^2(t)\ge\varepsilon^2$:

$$\frac{d}{dt}W(x(t)) \le -2c\left(x^2(t)-\varepsilon^2\right)^2 - 3c\varepsilon^2\left(x^2(t)-\varepsilon^2\right) + c^{-3}\|d\|_\infty^4 \tag{4.11}$$

Using (4.11) and definition (4.4) we get the following inequality for $t\in[r,t_{\max})$ a.e. when $x^2(t)\ge\varepsilon^2$:

$$\frac{d}{dt}W(x(t)) \le -4cW(x(t)) + c^{-3}\|d\|_\infty^4 \tag{4.12}$$

Using (4.4) and (4.5) we conclude that (4.12) holds even if $x^2(t)<\varepsilon^2$. Thus we conclude that (4.12) holds for $t\in[r,t_{\max})$ a.e..

Exploiting (4.12) and definition (4.4) we obtain the following estimate:

$$\left(x^2(t)-\varepsilon^2\right)^+ \le \exp\left(-2c(t-r)\right)\left(x^2(r)-\varepsilon^2\right)^+ + \frac{1}{\sqrt{2}c^2}\|d\|_\infty^2 \text{, for all } t\in[r,t_{\max}) \tag{4.13}$$

Estimate (4.13) gives the following estimate that holds when $t_{\max}>r$:

$$|x(t)| \le \exp\left(-c(t-r)\right)|x(r)| + \varepsilon + 2^{-1/4}c^{-1}\|d\|_\infty \text{, for all } t\in[r,t_{\max}) \tag{4.14}$$

Using (2.1), (2.3) and the fact that $p\left(x_t,u_t\right)\ge 0$, we also get for $t\in[0,t_{\max})$ a.e.:

$$x(t)\dot{x}(t) \le \left(\theta-2c\right)x^2(t) + x(t)d(t) \tag{4.15}$$

Using the inequality $x(t)d(t) \le \frac{1}{2\sqrt{2}c^2}d^2(t) + \frac{\sqrt{2}c^2}{2}x^2(t)$ in conjunction with (4.15) we get for $t\in[0,t_{\max})$ a.e.:

$$x(t)\dot{x}(t) \le \left(\theta + \frac{\sqrt{2}c^2}{2} - 2c\right)x^2(t) + \frac{1}{2\sqrt{2}c^2}\|d\|_\infty^2 \tag{4.16}$$

Differential inequality (4.16) gives the following estimate for all $t\in[0,t_{\max})$:

$$x^2(t) \le \exp\left(\left(2\theta+\sqrt{2}c^2-4c\right)^+ t\right)\left(x^2(0) + \frac{t}{\sqrt{2}c^2}\|d\|_\infty^2\right) \tag{4.17}$$

Combining (4.14) and (4.17) we get the following estimate when $t_{\max}>r$ that holds for all $t\in[r,t_{\max})$

$$|x(t)| \le M\exp\left(-c(t-r)^+\right)|x(0)| + \varepsilon + \left(1+\sqrt{r}M\exp\left(-c(t-r)^+\right)\right)2^{-1/4}c^{-1}\|d\|_\infty \tag{4.18}$$

where $M := \exp\left(\left(\theta + \frac{\sqrt{2}}{2}c^2 - 2c\right)^+ r\right)$. Taking into account estimate (4.17) we conclude that in every case ($t_{\max} > r$ or $t_{\max} \le r$) estimate (4.18) holds for all $t \in [0, t_{\max})$.

We next estimate the other component of the state, $u(t)$. Using (4.3), the fact that $\left(\theta \|x_t\|_2^2 + \langle x_t, d_t \rangle\right)^+ \le (\theta)^+ \|x_t\|_2^2 + \left|\langle x_t, d_t \rangle\right|$ and the Cauchy-Schwarz inequality $\left|\langle x_t, d_t \rangle\right| \le \|x_t\|_2 \|d_t\|_2 \le \sqrt{r}\, \|x_t\|_2 \|d_t\|_\infty$, we get for all $t \in [0, t_{\max})$:

$$|u(t)| \le 2c|x(t)| + \frac{2(\theta)^+ \|x_t\|_2^2 + 2\sqrt{r}\, \|x_t\|_2 \|d_t\|_\infty + cr|x(t)|^2}{2\left(\|x_t\|_2^2 + \left(\varepsilon^2 - x^2(t)\right)^+\right)} |x(t)| \tag{4.19}$$

Using (2.1), (2.3), definition (4.1) and the fact that $p(x_t, u_t) \ge 0$ (which implies that $x(t)u(t) \ge 0$ for $t \ge 0$), we also get for $t \in [-r, t_{\max})$ a.e.:

$$x(t)\dot{x}(t) \le (\theta)^+ x^2(t) + x(t)\left(d(t) + h(t)u(t)\right) \tag{4.20}$$

where

$$h(t) := \begin{cases} 1, & if\ t < 0 \\ 0, & if\ t \ge 0 \end{cases} \tag{4.21}$$

Using (4.20) we get for all $t, \tau \in [-r, t_{\max})$ with $t \ge \tau$:

$$x^2(t) \le \exp\left(2(\theta)^+ (t-\tau)\right) x^2(\tau) + 2\int_\tau^t \exp\left(2(\theta)^+ (t-s)\right) \left|x(s)\left(d(s) + h(s)u(s)\right)\right| ds \tag{4.22}$$

Employing (4.22), the fact that

$$\begin{aligned} &\int_\tau^t \exp\left(2(\theta)^+ (t-s)\right) \left|x(s)\left(d(s) + h(s)u(s)\right)\right| ds \\ &\le \exp\left(2(\theta)^+ (t-\tau)\right) \sup_{\tau \le s \le t}\left(\left|d(s) + h(s)u(s)\right|\right) \int_\tau^t |x(s)|\, ds \end{aligned}$$

and the Cauchy-Schwarz inequality $\int_\tau^t |x(s)|\, ds \le \sqrt{t-\tau} \left(\int_\tau^t |x(s)|^2 ds\right)^{1/2}$, we get for all $t, \tau \in [-r, t_{\max})$ with $t \ge \tau$:

$$\begin{aligned} x^2(t) &\le \exp\left(2(\theta)^+ (t-\tau)\right) x^2(\tau) \\ &+ 2\exp\left(2(\theta)^+ (t-\tau)\right) \sup_{\tau \le s \le t}\left(\left|d(s) + h(s)u(s)\right|\right) \sqrt{t-\tau} \left(\int_\tau^t |x(s)|^2 ds\right)^{1/2} \end{aligned} \tag{4.23}$$

Consequently, we obtain from (4.23) for all $t \in [0, t_{\max})$ and $\tau \in \left[t-r, t\right]$:

$$x^2(t) \le \exp\left(2\left(\theta\right)^+ r\right) x^2(\tau) + 2\exp\left(2\left(\theta\right)^+ r\right) \left\|(d+hu)_t\right\|_\infty \sqrt{r} \left\|x_t\right\|_2 \tag{4.24}$$

Integrating (4.24) with respect to $\tau \in \left[t-r, t\right]$, we get for all $t \in [0, t_{\max})$:

$$r x^2(t) \le \exp\left(2\left(\theta\right)^+ r\right) \left\|x_t\right\|_2^2 + 2\exp\left(2\left(\theta\right)^+ r\right) \left\|(d+hu)_t\right\|_\infty r\sqrt{r} \left\|x_t\right\|_2 \tag{4.25}$$

Using the inequality $2\sqrt{r} \exp\left(2\left(\theta\right)^+ r\right) \left\|(d+hu)_t\right\|_\infty \left\|x_t\right\|_2 \le b^{-1} r \exp\left(4\left(\theta\right)^+ r\right) \left\|(d+hu)_t\right\|_\infty^2 \left\|x_t\right\|_2^2 + b$ that holds for all $b > 0$, we get from (4.25) for all $t \in [0, t_{\max})$ and $b > 0$:

$$r x^2(t) \le \left(1 + b^{-1} r \left\|(d+hu)_t\right\|_\infty^2\right) r \exp\left(4\left(\theta\right)^+ r\right) \left\|x_t\right\|_2^2 + br \tag{4.26}$$

Combining (4.19) and (4.26) we get for all $t \in [0, t_{\max})$ and $b > 0$:

$$\begin{aligned} \left|u(t)\right| \le 2c\left|x(t)\right| &+ \frac{\left(2\left(\theta\right)^+ + c\left(1 + b^{-1} r \left\|(d+hu)_t\right\|_\infty^2\right) r \exp\left(4\left(\theta\right)^+ r\right)\right) \left\|x_t\right\|_2^2}{2\left(\left\|x_t\right\|_2^2 + \left(\varepsilon^2 - x^2(t)\right)^+\right)} \left|x(t)\right| \\ &+ \frac{2\sqrt{r} \left\|x_t\right\|_2 \left\|d_t\right\|_\infty + cbr}{2\left(\left\|x_t\right\|_2^2 + \left(\varepsilon^2 - x^2(t)\right)^+\right)} \left|x(t)\right| \end{aligned} \tag{4.27}$$

Using the inequality $2\sqrt{r} \left\|x_t\right\|_2 \left\|d_t\right\|_\infty \le c^{-1} b^{-1} \left\|d_t\right\|_\infty^2 \left\|x_t\right\|_2^2 + cbr$ we get from (4.27) for all $t \in [0, t_{\max})$ and $b > 0$:

$$\begin{aligned} \left|u(t)\right| \le 2c\left|x(t)\right| &+ \frac{cbr\left|x(t)\right|}{\left\|x_t\right\|_2^2 + \left(\varepsilon^2 - x^2(t)\right)^+} \\ &+ \frac{1}{2}\left(2\left(\theta\right)^+ + cr + b^{-1}\left(c^{-1}\left\|d_t\right\|_\infty^2 + cr^2 \left\|(d+hu)_t\right\|_\infty^2\right)\right) \exp\left(4\left(\theta\right)^+ r\right) \left|x(t)\right| \end{aligned} \tag{4.28}$$

Exploiting (3.1) and (4.28) we get for all $t \in [0, t_{\max})$ and $b > 0$:

$$\begin{aligned} \left|u(t)\right| &\le 2(3+r)\varepsilon^{-3}\left(r\sqrt{2}\left\|\dot{x}_t\right\|_\infty + \varepsilon\right) cb\left|x(t)\right| \\ &+ \frac{1}{2}\left(2\left(\theta\right)^+ + c(r+4) + b^{-1}\left(c^{-1}\left\|d_t\right\|_\infty^2 + cr^2 \left\|(d+hu)_t\right\|_\infty^2\right)\right) \exp\left(4\left(\theta\right)^+ r\right) \left|x(t)\right| \end{aligned} \tag{4.29}$$

Equation (2.1) implies that $\left\|\dot{x}_t\right\|_\infty \le \left|\theta\right| \left\|x_t\right\|_\infty + \left\|u_t\right\|_\infty + \left\|d_t\right\|_\infty$ for all $t \in [0, t_{\max})$. The previous inequality in conjunction with (4.29) gives for all $t \in [0, t_{\max})$ and $b > 0$:

$$
\begin{aligned}
&|u(t)| \le 2(3+r)\varepsilon^{-3} r\sqrt{2}\,\|u_t\|_\infty\, cb|x(t)| \\
&+\left((\theta)^+ + 2(3+r)\varepsilon^{-3}\left(r\sqrt{2}\left(|\theta|\|x_t\|_\infty + \|d_t\|_\infty\right)+\varepsilon\right)cb\right)\exp\left(4(\theta)^+ r\right)|x(t)| \\
&+\left(c(r+2)+c^{-1}b^{-1}\left(\|d_t\|_\infty^2 + c^2 r^2\,\|(d+hu)_t\|_\infty^2\right)\right)\exp\left(4(\theta)^+ r\right)|x(t)|
\end{aligned} \tag{4.30}
$$

Selecting $b = \dfrac{\varepsilon^3}{4\sqrt{2}(3+r)rc|x(t)|}$ when $|x(t)|>0$ we get from (4.30) for all $t\in[0,t_{\max})$:

$$
|u(t)| \le \frac{1}{2}\|u_t\|_\infty + \exp\left(4(\theta)^+ r\right) f(t) \tag{4.31}
$$

where

$$
\begin{aligned}
f(t) &:= \left(2|\theta| + c(r+3)\right)\|x_t\|_\infty + \frac{\|d_t\|_\infty}{2} + \frac{\varepsilon}{2\sqrt{2}r} \\
&+4\sqrt{2}(r+3)\varepsilon^{-3} r\left(\|d_t\|_\infty + cr\,\|(d+hu)_t\|_\infty\right)^2 |x(t)|^2
\end{aligned} \tag{4.32}
$$

Notice that (4.31) also holds when $|x(t)|=0$ (because in this case (4.30) implies that $|u(t)|=0$). Therefore, we conclude that (4.31) holds in any case for all $t\in[0,t_{\max})$.

Taking into account (4.21) and the fact that $u(s)=u_0(s)$ for all $s\in[-r,0]$, we conclude that $\|(d+hu)_t\|_\infty \le \|d_t\|_\infty + h(t-r)\|u_0\|_\infty$ for all $t\ge 0$. Taking into account (4.1) and the facts that $x(s)=x_0(s)$, $u(s)=u_0(s)$ for all $s\in[-r,0]$, we conclude that $\|d_t\|_\infty \le h(t-r)\left(\|u_0\|_\infty + \|\dot{x}_0\|_\infty + |\theta|\|x_0\|_\infty\right) + \|d\|_\infty$ for all $t\ge 0$. Since (4.18) gives for all $t\in[0,t_{\max})$

$$
\|x_t\|_\infty \le M\exp\left(-c(t-2r)^+\right)\|x_0\|_\infty + \varepsilon + \left(1+\sqrt{r}M\exp\left(-c(t-2r)^+\right)\right)2^{-1/4}c^{-1}\|d\|_\infty \tag{4.33}
$$

we obtain from (4.32) using the fact that $(a+b)^2 \le 2a^2+2b^2$ (repeatedly) the following inequality that holds for all $t\in[0,t_{\max})$:

$$
\begin{aligned}
f(t) &\le \left(2|\theta| + c(r+3)\right)M\exp\left(-c(t-2r)^+\right)\|x_0\|_\infty + \frac{h(t-r)}{2}\left(\|u_0\|_\infty + \|\dot{x}_0\|_\infty + |\theta|\|x_0\|_\infty\right) \\
&+K\left(\varepsilon^{-2}M^2\|x_0\|_\infty^2 + 2\right)h(t-r)\left(2\|u_0\|_\infty + \|\dot{x}_0\|_\infty + |\theta|\|x_0\|_\infty\right)^2 + \frac{\varepsilon}{2\sqrt{2}r} + \varepsilon\left(2|\theta| + c(r+3)\right) \\
&+\left(\left(2|\theta| + c(r+3)\right)\left(1+\sqrt{r}M\right)2^{-1/4}c^{-1} + \frac{1}{2}\right)\|d\|_\infty \\
&+K\left(2+\varepsilon^{-2}M^2\|x_0\|_\infty^2 + c^{-2}\varepsilon^{-2}\left(2\|u_0\|_\infty + \|\dot{x}_0\|_\infty + |\theta|\|x_0\|_\infty + \|d\|_\infty\right)^2\left(1+\sqrt{r}M\right)^2\right)\|d\|_\infty^2
\end{aligned} \tag{4.34}
$$

where $K = 16\varepsilon^{-1}\sqrt{2}(r+3)r(1+cr)^2$. Estimate (4.34) shows that the function $f(t)$ is bounded on $[0,t_{\max})$ with $f(t)\le F$ for all $t\in[0,t_{\max})$. Moreover, we obtain from (4.31) for all $t\in[0,t_{\max})$:

$$\max_{0\le s\le t}\left(|u(s)|\right)\le \frac{1}{2}\max_{-r\le s\le t}\left(|u(s)|\right)+\exp\left(4(\theta)^{+} r\right)F \tag{4.35}$$

Distinguishing the cases $\max_{-r\le s\le t}\left(|u(s)|\right)=\max_{0\le s\le t}\left(|u(s)|\right)$ and $\max_{-r\le s\le t}\left(|u(s)|\right)=\max_{-r\le s\le 0}\left(|u(s)|\right)$, using the fact that $u(s)=u_0(s)$ for all $s\in[-r,0]$, we get the following estimate for all $t\in[0,t_{\max})$:

$$\max_{0\le s\le t}\left(|u(s)|\right)\le \frac{1}{2}\|u_0\|_\infty+2\exp\left(4(\theta)^{+} r\right)F \tag{4.36}$$

Estimate (4.36) shows that $u(t)$ is bounded on $[0,t_{\max})$. It follows from Theorem 1 that $t_{\max}=+\infty$. Furthermore, estimate (2.11) is a direct consequence of estimate (4.18).

We next finish the proof by showing the validity of estimate (2.12). Let $\omega\in(0,c]$ be such that $\exp(\omega r)<2$. Using (4.31) we get for all $t\ge 0$:

$$|u(t)|\exp(\omega t)\le \frac{\exp(\omega r)}{2}\max_{t-r\le s\le t}\left(|u(s)|\exp(\omega s)\right)+\exp\left(4(\theta)^{+} r\right)f(t)\exp(\omega t) \tag{4.37}$$

Inequality (4.37) implies the following estimate for all $t\ge 0$:

$$\begin{aligned}\max_{0\le s\le t}\left(|u(s)|\exp(\omega s)\right)&\le \frac{\exp(\omega r)}{2}\max_{-r\le s\le t}\left(|u(s)|\exp(\omega s)\right)\\ &+\exp\left(4(\theta)^{+} r\right)\max_{0\le s\le t}\left(f(s)\exp(\omega s)\right)\end{aligned} \tag{4.38}$$

By distinguishing the cases $\max_{-r\le s\le t}\left(|u(s)|\exp(\omega s)\right)=\max_{0\le s\le t}\left(|u(s)|\exp(\omega s)\right)$ and $\max_{-r\le s\le t}\left(|u(s)|\exp(\omega s)\right)=\max_{-r\le s\le 0}\left(|u(s)|\exp(\omega s)\right)$, we obtain from (4.38) for all $t\ge 0$:

$$\begin{aligned}\max_{0\le s\le t}\left(|u(s)|\exp(\omega s)\right)&\le \frac{\exp(\omega r)}{2}\max_{-r\le s\le 0}\left(|u(s)|\exp(\omega s)\right)\\ &+2\left(2-\exp(\omega r)\right)^{-1}\exp\left(4(\theta)^{+} r\right)\max_{0\le s\le t}\left(f(s)\exp(\omega s)\right)\end{aligned} \tag{4.39}$$

Thus, we obtain from (4.39) and the fact that $u(s)=u_0(s)$ for all $s\in[-r,0]$, the following inequality for all $t\ge 0$:

$$|u(t)|\le \frac{\exp(-\omega(t-r))}{2}\|u_0\|_\infty+2\left(2-\exp(\omega r)\right)^{-1}\exp\left(4(\theta)^{+} r\right)\max_{0\le s\le t}\left(f(s)\exp(-\omega(t-s))\right) \tag{4.40}$$

Using the facts that $\omega\in(0,c]$, $\exp\left(-c(t-2r)^{+}\right)\le\exp\left(-c(t-2r)\right)$, $h(t-r)\le\exp\left(-c(t-r)\right)$ in conjunction with (4.34) we get:

$$
\begin{aligned}
&\max_{0\le s\le t}\left( f(s)\exp(-\omega(t-s))\right) \le \frac{\varepsilon}{2\sqrt{2r}} + \varepsilon\left(2|\theta| + c(r+3)\right) \\
&+\left(3|\theta| + c(r+3)\right) M \exp(-\omega t)\exp\left(2cr\right)\left( \|x_0\|_\infty + \|u_0\|_\infty + \|\dot{x}_0\|_\infty \right) \\
&+K\left(\varepsilon^{-2}M^2\|x_0\|_\infty^2 + 2\right)\left(2+|\theta|\right)^2 \exp(-\omega t)\exp\left(cr\right)\left(\|u_0\|_\infty + \|\dot{x}_0\|_\infty + \|x_0\|_\infty\right)^2 \\
&+\left(\left(2|\theta| + c(r+3)\right)\left(1+\sqrt{r}M\right)2^{-1/4}c^{-1} + \frac{1}{2}\right)\|d\|_\infty \\
&+K\left(2+\varepsilon^{-2}\left(M^2 + c^{-2}\left(2+|\theta|\right)^2\left(1+\sqrt{r}M\right)^2\right)\left(\|u_0\|_\infty + \|\dot{x}_0\|_\infty + \|x_0\|_\infty + \|d\|_\infty\right)^2\right)\|d\|_\infty^2
\end{aligned}
\tag{4.41}
$$

Estimate (2.12) for certain appropriate function $a \in K_\infty$ and certain non-decreasing function $\gamma : \mathbb{R}_+ \to \mathbb{R}_+$ is a direct consequence of estimates (4.40), (4.41).

The proof is complete. $\lhd$

We end this section by providing the proof of Theorem 3.

**Proof of Theorem 3:** Let $\left(x_0, u_0, \hat{\theta}_0\right) \in X \times \mathbb{R}$, $d \in L^\infty\left(\mathbb{R}_+\right)$ be given. By virtue of Theorem 2 there exists a unique pair of continuous mappings $x, u : \mathbb{R}_+ \to \mathbb{R}$ with $x \in W^{1,\infty}\left(\mathbb{R}_+\right)$, $x(s) = x_0(s)$, $u(s) = u_0(s)$ for all $s \in [-r, 0]$ for which (2.1) holds for $t \ge 0$ a.e., (2.3) holds for all $t \ge 0$. It follows from (2.4), (2.5), (2.6) that there exists a unique piecewise constant function $\hat{\theta} : \mathbb{R}_+ \to \mathbb{R}$ with $\hat{\theta}(0) = \hat{\theta}_0$ which satisfies (2.4) for all $t \in \left[ir, (i+1)r\right)$, $i \in \mathbb{N}$ and (2.5), (2.6) for all $i \in \mathbb{N}$.

We define the time sets:

$$
\begin{aligned}
I &= \left\{ s \ge r : \|x_s\|_2 \ge \sigma \right\} \\
J &= \left\{ s \in \left[0, r\right) : \|x_s\|_2 \ge \sigma \right\}
\end{aligned}
\tag{4.42}
$$

Due to (2.8), (4.2), (4.42), (2.4), (2.5), (2.6), we obtain in the case $I \ne \varnothing$:

$$
\left|\hat{\theta}(t) - \theta\right| \le \frac{\sqrt{r}}{\sigma}\|d\|_\infty \text{ for all } t \ge r\left\lceil r^{-1}\min(I)\right\rceil \tag{4.43}
$$

$$
\begin{gathered}
\hat{\theta}(t) = \hat{\theta}(0) \text{ for } 0 \le t < r\left\lceil r^{-1}\min(I)\right\rceil, \\
\text{when } J = \varnothing \text{ or } J \ne \varnothing \text{ and } \left\lceil r^{-1}\min(I)\right\rceil = 1
\end{gathered}
\tag{4.44}
$$

$$
\begin{gathered}
\hat{\theta}(t) = \hat{\theta}(0) \text{ for } 0 \le t < r \text{ and} \\
\hat{\theta}(t) = q\left(x_\tau, u_\tau\right) \text{ for } r \le t < r\left\lceil r^{-1}\min(I)\right\rceil \text{ with } \tau = \max\left\{ t \in [0, r] : \|x_t\|_2 = \max_{0\le s\le r}\left(\|x_s\|_2\right) \right\}, \\
\text{when } J \ne \varnothing \text{ and } \left\lceil r^{-1}\min(I)\right\rceil \ge 2
\end{gathered}
\tag{4.45}
$$

Similarly, due to (4.42), (2.4), (2.5), (2.6), we obtain in the case $I = \varnothing$:

$$\hat{\theta}(t) = \hat{\theta}(0) \text{ for } t \ge 0 \text{, when } J = \varnothing \tag{4.46}$$

$$\hat{\theta}(t) = \hat{\theta}(0) \text{ for } 0 \le t < r \text{ and}$$
$$\hat{\theta}(t) = q\left(x_\tau, u_\tau\right) \text{ for } t \ge r \text{ with } \tau = \max\left\{ t \in [0,r] : \|x_t\|_2 = \max_{0 \le s \le r}\left(\|x_s\|_2\right) \right\},$$
$$\text{when } J \ne \varnothing \tag{4.47}$$

Therefore, we conclude that $\hat{\theta} \in L^\infty\left(\mathbb{R}_+\right)$. Estimate (2.13) with $T = r\left\lceil r^{-1}\min(I)\right\rceil$ in the case where there exists $t \ge r$ with $\|x_t\|_2 \ge \sigma$ (i.e., in the case $I \ne \varnothing$; recall definition (4.42)) is a direct consequence of (4.43). The proof is complete. ⊲

## 5. Concluding Remarks

As noticed earlier, this paper is the first step towards the resolution of several decades-old challenges in disturbance-robust adaptive control. We have showed in a scalar linear system that the use of delays in adaptive control can guarantee features that cannot be guaranteed by adaptive controllers without delays. However, many things are left to be done in the future. First, the results have to be extended to important classes of systems (e.g., systems in strict feedback form or systems satisfying the matching condition; see [16,17,20]). This extension is far from trivial and will require the development of novel design procedures as well as the generation of new mathematical results that will make the analysis of the resulting closed-loop systems tractable. Second, simulation studies have to be performed to investigate how the performance is affected by different choices for the controller parameters. Indeed, the obtained stability estimates show us some important qualitative characteristics but the actual dependence of the performance on the controller parameters can be fully understood only with the combination of theory and various numerical experiments. The numerical studies should also compare the proposed delay-adaptive scheme with well-established existing schemes (e.g., leakage).